\documentclass{article}

\usepackage{amsmath}	
\usepackage{amsfonts}
\usepackage{amssymb}
\usepackage{booktabs}
\usepackage{hyperref}
\usepackage[capitalize]{cleveref}
\crefname{equation}{}{}
\usepackage{tikz}
\usetikzlibrary{patterns, spy, calc}
\usepackage{pgfplots}
\pgfplotsset{compat=1.4}
\pgfdeclarelayer{back}
\pgfsetlayers{back,main}

\tikzset{%
  on layer/.code={
    \pgfonlayer{#1}\begingroup
    \aftergroup\endpgfonlayer
    \aftergroup\endgroup
  },%
  draw on back/.style={
    preaction={
      draw,
      on layer=back,
      thick
    },
    gray!30
  }
}

\newcommand{\ooc}{1\!/\!c}

\title{Localized Reduced Basis Approximation of a Nonlinear 
   Finite Volume Battery Model with Resolved Electrode Geometry\thanks{This work has been supported by the German Federal Ministry of Education and Research (BMBF) under contract number 05M13PMA.}}
\author{Mario Ohlberger\thanks{Mario Ohlberger and Stephan Rave, Applied Mathematics Muenster \& Center for Nonlinear
	Science, University of Muenster, Einsteinstr. 62, 48149 Muenster, Germany.  Email:
\texttt{(mario.ohlberger|stephan.rave)@uni-muenster.de}} \and Stephan Rave\footnotemark[1]}
\date{June 15, 2016}

\begin{document}
\maketitle

\begin{abstract}
	In this contribution we present first results towards localized model
	order reduction for spatially resolved, three-dimensional lithium-ion
	battery models.
	We introduce a localized reduced basis scheme based on non-conforming 
	local approximation spaces stemming from a finite volume discretization
	of the analytical model and localized empirical operator interpolation
	for the approximation of the model's nonlinearities.
	Numerical examples are provided indicating the feasibility of our
	approach.
\end{abstract}

\section{Introduction}
\label{sec:introduction}
Over the recent years, three dimensional lithium (Li) ion battery models that fully
resolve the microscopic geometry of the battery electrodes have become a
subject of active research in electrochemistry \cite{LatzZausch2011}.
These models are also studied in the collaborative research project \textsc{Multibat},
where the influence of the microscopic electrode geometry plays in
important role in understanding the degradation process of Li-plating \cite{HeinLatz2016}.

Due to the strongly nonlinear character of these models and the large number of
degrees of freedom of their discretization, numerical simulations are time consuming
and parameter studies quickly turn prohibitively expensive.
Our work in context of the \textsc{Multibat} project has shown that
model reduction techniques such as reduced
basis (RB) methods are able to vastly reduce the computational complexity of parametrized microscale
battery models while retaining the full microscale features of their solutions
\cite{OhlbergerRaveEtAl2016,OhlbergerRaveEtAl2014}.
Still, such methods depend on the solution of the full high-dimensional model for selected
parameters during the so-called offline phase.
When only relatively few simulations of the model are required -- as it is typically the case
for electrochemistry simulations where one is mainly interested in the qualitative behaviour
of the battery cell -- the offline phase can quickly take nearly as much time as the simulation
of the full model for all parameters of interest would have required.
It is, therefore, paramount to reduce the number of full model solves as much possible.
Localized RB methods construct spatially localized approximation spaces from
few global model solves or even by only solving adequate local problems
(see also \cite{BuhrEngwerEtAl2015,BuhrEngwerEtAl2014,OhlbergerRaveEtAl2015} and the references therein). 
Thus, these methods are a natural choice for accelerating the offline phase of RB schemes,
in particular for problems with a strong microscale character such as geometrically resolved
electrochemistry simulations.

While localized RB methods have been studied extensively for linear problems and while there are
first results for instationary problems \cite{OhlbergerRaveEtAl2015,OhlbergerRaveEtAl2016}, we are not
aware of any previous work treating nonlinear models.
In this contribution, we introduce a localized RB scheme for nonlinear finite volume battery
models, which builds local approximation and interpolation spaces by decomposition of global
solution snapshots w.r.t.\ a given coarse triangulation of the computational domain
(\cref{sec:lrb}).
As a preparation, we will first briefly review the microscale battery model under
consideration (\cref{sec:model}), its discretization (\cref{sec:discretization}) and finally its
RB approximation (\cref{sec:rb}).
We will close with first numerical experiments that investigate the applicability of localized RB techniques
to the problem at hand (\cref{sec:experiment}).

\section{Analytical model}\label{sec:model}

\newcommand\irregularellipse[3]{
\pgfextra {\pgfmathsetmacro{\lena}{(#1)}\pgfmathsetmacro{\lenb}{\lena*(#2)}}
+(0:{\lena} and \lenb)
  \foreach \a in {10,20,...,350}{
  \pgfextra {\pgfmathsetmacro{\lena}{(#1) +rand*(#3)}\pgfmathsetmacro{\lenb}{\lena*(#2)}}
  -- +(\a:{\lena} and \lenb)
  } -- cycle
}

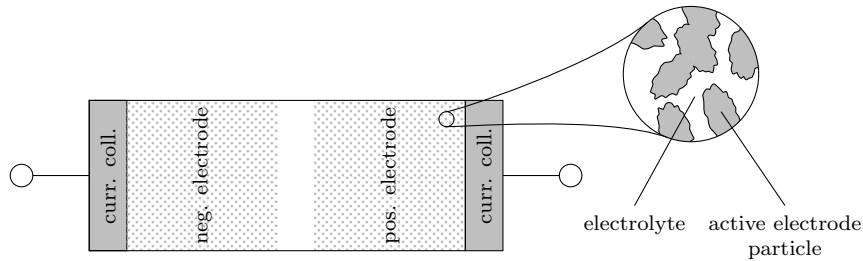
\begin{figure}[t]
	\centering
	\begin{tikzpicture}[scale=0.5, font=\footnotesize]
		\draw[fill=lightgray] (0,0) rectangle (1, 4);	
		\node[rotate=90] at (0.5, 2) {curr. coll.};
		\fill[pattern=crosshatch dots, pattern color=lightgray](1,0) rectangle (5, 4);	
		\node[rotate=90] at (3, 2) {neg. electrode};
		\node[rotate=90] at (6, 2) {};
		\fill[pattern=crosshatch dots, pattern color=lightgray] (6,0) rectangle (10, 4);	
		\node[rotate=90] at (8, 2) {pos. electrode};
		\draw[fill=lightgray] (10,0) rectangle (11, 4);	
		\node[rotate=90] at (10.5, 2) {curr. coll.};
		\draw (0,0) rectangle (10, 4);	
		\draw (0,2) -- (-1.5, 2);
		 \draw (-1.8, 2) circle (0.3);
		\draw (11,2) -- (12.5, 2);
		\draw (12.8, 2) circle (0.3);
		\coordinate (A) at (16, 4.7);
		\pgfmathsetmacro{\Aradius}{1.8}
		\coordinate (Alow) at ($(A) + (-110:\Aradius)$);
		\coordinate (Aup) at ($(A) + (130:\Aradius)$);
		\coordinate (B) at (9.5, 3.5);
		\pgfmathsetmacro{\Bradius}{0.2}
		\coordinate (Blow) at ($(B) + (-100:\Bradius)$);
		\coordinate (Bup) at ($(B) + (90:\Bradius)$);
		\begin{scope}
			\clip (A) circle (\Aradius);
			\clip[on layer=back] (A) circle (\Aradius);
			\fill[draw on back, fill=lightgray, rotate=69, rounded corners=0.3pt] ($(A) + (.6, 1.7)$) \irregularellipse{0.4}{1.8}{0.05};
			\fill[draw on back, fill=lightgray, rotate=24, rounded corners=0.3pt] ($(A) + (0.9, 1.4)$) \irregularellipse{0.4}{1.8}{0.05};
			\fill[draw on back, fill=lightgray, rotate=34, rounded corners=0.3pt] ($(A) + (0.5, 0.5)$) \irregularellipse{0.4}{1.8}{0.05};
			\fill[draw on back, fill=lightgray, rotate=-4, rounded corners=0.3pt] ($(A) + (1.4, 0.7)$) \irregularellipse{0.4}{1.8}{0.05};
			\fill[draw on back, fill=lightgray, rotate=-44, rounded corners=0.3pt] ($(A) + (-0.3, -0.4)$) \irregularellipse{0.38}{1.8}{0.05};
			\fill[draw on back, fill=lightgray, rotate=20, rounded corners=0.3pt] ($(A) + (0.4, -1.25)$) \irregularellipse{0.4}{1.8}{0.05};
			\fill[draw on back, fill=lightgray, rotate=20, rounded corners=0.3pt] ($(A) + (-0.9, -1.25)$) \irregularellipse{0.4}{1.8}{0.05};
		\end{scope}
		\draw (A) circle (\Aradius);
		\draw (B) circle (\Bradius);
		\draw (Alow) .. controls ($(Alow)!1!90:(A)$) and ($(Blow)!1!-90:(B)$) .. (Blow);
		\draw (Aup) .. controls ($(Aup)!1!-90:(A)$) and ($(Bup)!1!90:(B)$) .. (Bup);
		\draw ($(A) + (0.8, -1.0)$) -- ($(A) + (2.5, -3.5)$) coordinate [label={[align=center]below:active electrode\\particle}];
		\draw ($(A) + (0.1, -0.6)$) -- ($(A) + (-1.5, -3.5)$) coordinate [label=below:electrolyte];
	\end{tikzpicture}
	\caption{Sketch of a lithium-ion battery cell. The cell is connected via two metallic current collectors which are in
	contact with the negative/positive cell electrodes. The porous electrodes are composed of active electrode
	particles into which Li-ions intercalate from the electrolyte filling the pore space of the electrodes. }
	\label{fig:battery_schematic}
\end{figure}

As in \cite{OhlbergerRaveEtAl2016,OhlbergerRaveEtAl2014}, we consider the microscale 
battery model introduced in \cite{LatzZausch2011} (without taking thermal effects into account
and assuming constant $t_+$).
In this model, the battery cell is described via coupled systems of partial differential
equations for the concentration of Li$^+$-ions and the electrical potential $\phi$ for
each part of the cell: the electrolyte, positive and negative electrode, as well as
positive and negative current collector.

In the electrolyte, the change of the concentration $c$ is governed by a diffusion
process, whereas $\phi$ is determined by a stationary potential equation with source term
depending non-linearly on $c$:
\begin{align*}
\frac{\partial c}{\partial t} - \nabla \cdot (D_e \nabla c) &= 0,
\\
 -\nabla \cdot \Bigl(\kappa \frac{1 - t_+}{F}RT \frac{1}{c} \nabla c - \kappa \nabla
 \phi\Bigr) &= 0.
\end{align*}
In the electrodes, the evolution of $c$, i.e. the intercalation of Li-ions into the active
particles, is again driven by diffusion. The potential $\phi$ no longer depends
on the Li-ion distribution:
\begin{align*}
\frac{\partial c}{\partial t} - \nabla \cdot (D_s \nabla c) &= 0,
\\
-\nabla \cdot (\sigma \nabla \phi) &= 0.
\end{align*}
No Li-ions enter the metallic current collectors, so $c \equiv 0$ on this part of the
domain, whereas $\phi$ is again given as:
\begin{equation*}
	-\nabla \cdot (\sigma \nabla \phi) = 0.
\end{equation*}

The reaction at the interface between active electrode particles and the electrolyte
is governed by the so-called Butler-Volmer kinetics which determine the electric
current $j = \nabla \phi \cdot n$ from the active particle into the electrolyte 
as
\begin{equation}\label{eq:butler_volmer}
       j = 2k\sqrt{c_ec_s(c_{max}- c_s)} \sinh\left(\frac{\phi_s - \phi_e -
		       U_0(\frac{c_s}{c_{max}})}{2RT} \cdot F \right),
\end{equation}
where $c_e$, $\phi_e$ ($c_s$, $\phi_s$) are the concentration and potential on the
electrolyte (solid particle) side of the interface.
The Li-ion flux $N$ over the interface proportionally depends on $j$ and is given
by $N = j / F$.
Note that the Butler-Volmer relations ensure the coupling between both considered
variables and, through the exponential dependence on the overpotential
$\phi_s - \phi_e - U_0(c_s/c_{max})$, lead to a highly nonlinear behaviour of the system.

Finally, continuity conditions for $\phi$ are imposed between electrode and current
collector, whereas there is no coupling between electrolyte and current collector.
The following boundary conditions are imposed: $\phi = U_0(c(0) / c_{max})$ at the negative
current collector boundary, Neumann boundary conditions at the positive current
collector (applied fixed charge/discharge rate) and periodic boundary conditions
for $c$ and $\phi$ at the remaining domain boundaries.
We denote the initial concentration at time $t=0$ by $c_0 = c(0)$, the final time
is denoted as $T$.
All appearing natural/material constants as well as the initial data is summarised
in \cref{tab:constants}.

\section{Finite volume discretization}\label{sec:discretization}

Following \cite{PopovVutovEtAl2011}, the continuous model is discretized using a basic cell centered finite
volume scheme on a voxel grid. Each voxel is assigned a
unique subdomain and the Butler-Volmer conditions are chosen as numerical flux on grid faces
separating an electrolyte from an electrode voxel.
We obtain a single nonlinear finite volume operator $A_\mu: V_h \oplus V_h \to V_h \oplus V_h$
for the whole computational domain, where $V_h$ denotes the space of piecewise constant grid functions
and $\mu$ indicates a parameter we want to vary.
In the following, we will consider the applied charge current as parameter of interest.
Implicit Euler time stepping with constant time step size $\Delta t$ leads to the $N:=T/\Delta t$
nonlinear equation systems
\begin{equation}
\label{eq:detailed}
        \begin{bmatrix}
                 \frac{1}{\Delta t}(c_{\mu}^{(n+1)} - c_{\mu}^{(n)}) \\
                 0
        \end{bmatrix}
         + A_\mu
         \left(\begin{bmatrix}
                c_{\mu}^{(n+1)} \\
                \phi_{\mu}^{(n+1)}
         \end{bmatrix}\right)
         = 0, \qquad (c_{\mu}^{(n)}, \phi_{\mu}^{(n)}) \in V_h \oplus V_h.
\end{equation}
The equation systems are solved using a standard Newton iteration scheme.

\begin{table}
\caption{Constants used in the battery model. The open circuit potential $U_0$
for a state of charge $s$ is give as $U_0(s) = (-0.132+1.41\cdot e^{-3.52 s})\cdot V$
for the negative electrode and as $U_0(s) = [0.0677504 \cdot \tanh (-21.8502\cdot s+ 12.8268)
- 0.105734 \cdot \bigl( (1.00167 - s)^{-0.379571} - 1.576\bigr) - 0.045\cdot e^{-71.69 \cdot s^8}
+ 0.01\cdot e^{-200\cdot(s-0.19)} + 4.06279]\cdot V$ for the positive electrode.}\label{tab:constants}
\renewcommand{\arraystretch}{1}\setlength{\tabcolsep}{0.5em}
{\footnotesize
	\begin{tabular}{llll}
		\toprule
		symbol & unit & value & description \\
		\midrule
		$D_e$ & $\frac{cm^2}{s}$ & $1.622\cdot10^{-6}$ & collective interdiffusion coefficient in electrolyte \\[\medskipamount]
		$D_s$ & $\frac{cm^2}{s}$ & $10^{-10}$ & ion diffusion coefficient in electrodes \\[\medskipamount]
		$\sigma$ & $\frac{s}{cm}$ & $10$  & electronic conductivity in neg.\ electrode \\
		& & $0.38 $ & electronic conductivity in pos.\ electrode \\[\medskipamount]
		$\kappa$ & $\frac{s}{cm}$ & $0.02$ & ion conductivity \\[\medskipamount]
		$c_{max}$ & $\frac{mol}{cm^3}$ & $24681\cdot 10^{-6}$ & maximum Li$^+$ concentration in neg.\ electrode \\
		& & $23671 \cdot 10^{-6}$ & maximum Li$^+$ concentration in pos.\ electrode \\[\medskipamount]
		$k$ & $\frac{A cm^{2.5}}{mol^{1.5}}$ & $0.002 $ & reaction rate at neg. electrode/electrolyte interface \\
		& & $0.2$ & reaction rate at pos. electrode/electrolyte interface \\[\medskipamount]
		$c_0$ & $\frac{mol}{cm^3}$ & $1200 \cdot 10^{-6} $ & initial concentration in electrolyte \\
		& & $2639 \cdot 10^{-6}$ & initial concentration in neg.\ electrode \\
		& & $20574 \cdot 10^{-6}$ & initial concentration in pos.\ electrode \\[\medskipamount]
		$t_+$ & & $0.39989$ & transference number \\[\medskipamount]
		$T$ & $K$ & $298$ & temperature \\[\medskipamount]
		$F$ & $ \frac{As}{mol}$ & $96487$ & Faraday constant \\[\medskipamount]
		$R$ & $ \frac{J}{mol\,K}$ & $8.314$ & universal gas constant \\[\medskipamount]
		\bottomrule
	\end{tabular}
}
\end{table}

Note that we can decompose $A_\mu$ as
\begin{equation}\label{eq:operator_decomposition}
A_\mu = A_\mu^{(aff)}  + A^{(bv)} + A^{(\ooc)} 
\end{equation}
where $A^{(bv)}, A^{(\ooc)}: V_h \oplus V_h \to V_h \oplus V_h$ accumulate the numerical fluxes
corresponding to \cref{eq:butler_volmer} and $\kappa \frac{1 - t_+}{F}RT \frac{1}{c} \nabla c$.
Thus, the operator $A_\mu^{(aff)}$ collecting the remaining numerical fluxes is 
affine linear and the only operator in the decomposition depending on the charge rate.
$A_\mu^{(aff)}$ can be further decomposed as
\begin{equation}\label{eq:affine_operator_decomposition}
	A_\mu^{aff} = A^{const} + \mu\cdot A^{bnd} + A^{lin},
\end{equation}
with constant, non-parametric operators $A^{const}, A^{bnd}$ corresponding to the boundary
conditions and a non-parametric linear operator $A^{(lin)}$.

\section{Reduced basis approximation}\label{sec:rb}

As reduced model we consider the Galerkin projection of \cref{eq:detailed} onto an
appropriate reduced basis space $\tilde{V} \subseteq V_h \oplus V_h$, i.e.
we solve
\begin{equation}
\label{eq:reduced}
		P_{\tilde{V}} \left\{
                \begin{bmatrix}
                         \frac{1}{\Delta t}(\tilde{c}_{\mu}^{(n+1)} - \tilde{c}_{\mu}^{(n)}) \\
                         0
                \end{bmatrix} 
                 + 
                 A_\mu
                 \left(\begin{bmatrix}
		        \tilde{c}_{\mu}^{(n+1)} \\
			\tilde{\phi}_{\mu}^{(n+1)}
                 \end{bmatrix}\right) \right\}
         = 0, \quad (\tilde{c}_{\mu}^{(n)}, \tilde{\phi}_{\mu}^{(n)}) \in
	 \tilde{V},
\end{equation}
where $P_{\tilde{V}}$ denotes the orthogonal projection onto $\tilde{V}$.
In order to obtain at an online efficient scheme, the projected operator 
$P_{\tilde{V}} \circ A_\mu$ has to be approximated by an efficiently
computable approximation.
Considering the decompositions \cref{eq:operator_decomposition,eq:affine_operator_decomposition}, 
only the nonlinear operators $A_\mu^{(bv)}$, $A_\mu^{(\ooc)}$ require special treatment
for which we employ empirical operator interpolation \cite{HaasdonkOhlbergerEtAl2008}
based on the empirical interpolation method \cite{BarraultMadayEtAl2004}.
Denoting the discrete time differential operator by $B$, the fully reduced scheme
is then given as
\begin{equation}\label{eq:ei_reduced}
\begin{aligned}
	\Bigl\{P_{\tilde{V}} \circ B & +
		 P_{\tilde{V}} \circ A^{(const)} + \mu \cdot P_{\tilde{V}}
		 \circ A^{(bnd)} 
		  + P_{\tilde{V}} \circ A^{(lin)}  \\
		 & + \{P_{\tilde{V}} \circ I_{M^{(\ooc)}}^{(\ooc)}\} \circ
		 \tilde{A}_{M^{(\ooc)}}^{(\ooc)} \circ R_{M^{\prime (\ooc)}}^{(\ooc)} \\
		 & + \{P_{\tilde{V}} \circ I_{M^{(bv)}}^{(bv)}\} \circ
		 \tilde{A}_{M^{(bv)}}^{(bv)} \circ R_{M^{\prime (bv)}}^{(bv)}\  \ \  \Bigr\}
                 \left(\begin{bmatrix}
		        \tilde{c}_{\mu}^{(t+1)} \\
			\tilde{\phi}_{\mu}^{(t+1)}
                 \end{bmatrix}\right)
         = 0, 
 \end{aligned}
\end{equation}
where $\tilde{A}^{(*)}_{M^{(*)}}: \mathbb{R}^{M^{\prime(*)}} \to \mathbb{R}^{M^{(*)}}$ ($* = bv,\ooc$)
denotes the restriction of $A^{(*)}$ to certain $M^{(*)}$ image degrees of freedom given the 
required $M^{\prime(*)}$ source degrees of freedom, $R^{(*)}_{M^{\prime(*)}}: V_h \oplus V_h \to \mathbb{R}^{M^{\prime(*)}}$
is the linear operator restricting finite volume functions to these $M^{\prime(*)}$ source degrees of
freedom, and $I^{(*)}_{M^{(*)}}: \mathbb{R}^{M^{(*)}} \to V_h \oplus V_h$ is the linear interpolation
operator to the $M^{(*)}$ evaluation points and an appropriately selected interpolation basis.
Note that for the considered finite volume scheme we have $M^{\prime(*)} \leq 14\cdot M^{(*)}$,
thus $\tilde{A}^{(*)}_{M^{(*)}}$ can be computed quickly for sufficiently small $M^{(*)}$.
The remaining terms in \cref{eq:ei_reduced} are linear (or constant) and can be pre-computed for
a given reduced basis of $\tilde{V}$.

\begin{figure}[t]
	\begin{center}
		\raisebox{0.0\height}{\includegraphics[trim=0 80 80 20, clip,
		width=0.60\textwidth]{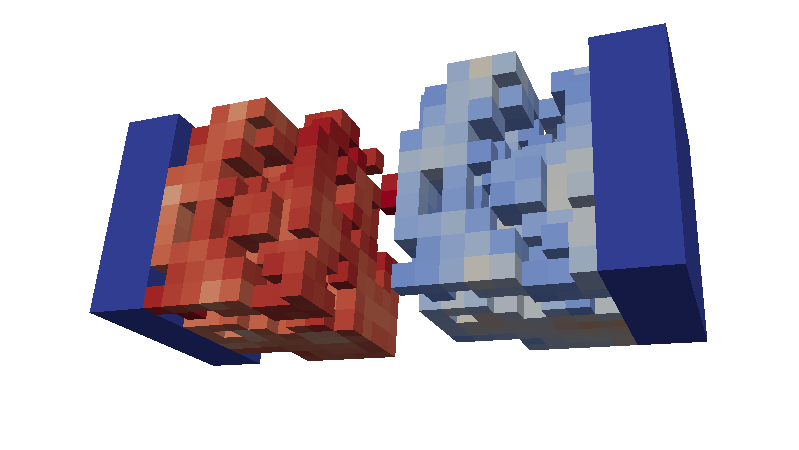}}
	\end{center}
\caption{
Small battery geometry used in numerical experiment.
Domain: $104\mu m \times 40\mu m \times 40\mu m$, $4.600$ degrees of freedom.
Coloring: Li$^+$-concentration at final simulation time $T=2000s$, electrolyte not depicted.
}\label{fig:geometry}
\end{figure}

In \cite{OhlbergerRaveEtAl2016} we have considered the solution of \cref{eq:ei_reduced} where $\tilde{V}$
and the interpolation data for $A^{(bv)}$, $A^{(\ooc)}$ have been generated using standard
model order reduction techniques.
The reduced space $\tilde{V}$ was determined by computing a proper orthogonal decomposition (POD) \cite{sirovich87}
of solution trajectories of \cref{eq:detailed} for an equidistant training set of charge
rate parameters.
Since the $c$ and $\phi$ variables are defined on different scales, the POD had to be applied
separately for both variables, yielding a reduced space of the form $\tilde{V} = \tilde{V}_c \oplus \tilde{V}_\phi$,
in order to obtain a stable scheme.
Moreover, the intermediate stages of the Newton algorithms used for solving \cref{eq:detailed}
were included in the snapshot data to ensure the convergence of the Newton algorithms when
solving the reduced scheme.

The interpolation bases and interpolation points have been obtained by evaluating $A^{(*)}$ on the
computed solution trajectories and then performing the \textsc{EI-Greedy} algorithm \cite{DrohmannHaasdonkEtAl2012}
on these evaluations.
Note that for solution trajectories of \cref{eq:detailed}, $A_\mu$ vanishes identically in the
$\phi$-component. Thus, applying the \textsc{EI-Greedy} algorithm directly to evaluations of $A_\mu$
would not have yielded usable interpolation spaces.

\section{Localized basis generation}\label{sec:lrb}

Localized RB methods can be seen as RB schemes where the reduced space
$\tilde{V}$ has a certain direct sum decomposition $\tilde{V} = \tilde{V}_1 \oplus \ldots \oplus \tilde{V}_K$
with subspaces $\tilde{V}_i$ associated to some partition $\overline{\Omega} = 
\overline{\Omega_1} \cup \ldots \cup \overline{\Omega_K}$ of the computational domain $\Omega$.
Since this imposes an additional constraint on the possible choices of reduced spaces $\tilde{V}$,
it is not to be expected that such methods yield better approximation spaces for the same (total)
dimension of $\tilde{V}$ than classical RB methods.
However, these methods can yield enormous saving in computation time during basis generation.
In particular, we are interested in the following aspects:
\begin{enumerate}
	\item When the parametrization of the problem mainly affects the global solution behaviour,
        	only few global solution snapshots may be required to observe all relevant local
		behaviour.
		This can be exploited by computing local approximation spaces from global
		solutions which have been decomposed according to the partition
		$\Omega_1 \cup \ldots \cup \Omega_K$ (e.g.~\cite{AlbrechtHaasdonkEtAl2012}).
	\item The local approximation spaces $\tilde{V}_i$ may be enriched by solving appropriate
		local problems on a neighbourhood of $\Omega_i$.
		The solution of the local problems can be trivially parallelized, and each local
		problem will be solvable much faster than the global problem, which might even
		be unsolvable with the available computational resources (e.g.~\cite{OhlbergerSchindler2015}).
	\item When the problem undergoes local changes (e.g. geometry change due to Li-plating),
		the spaces $\tilde{V}_i$ which are unaffected
		by the change can be reused and only few new local problems have to be solved
		(e.g.~\cite{BuhrEngwerEtAl2015}).
\end{enumerate}

For many applications, the time for basis generation must be taken into account when considering
the overall efficiency of the reduction scheme. Hence, such localization 
approaches can be an essential tool for making model order reduction profitable for these
applications.
This is also the case for battery simulations, where typically only relatively few parameter
samples are required to gain an appropriate idea of the behaviour of the model and these same
computational resources are available for all required simulations.
Also note that while reduced system matrices/Jacobians are dense matrices for standard RB
approaches, one typically obtains block sparse matrices for localized RB approaches,
so the increased global system dimension can be largely compensated by appropriate choices
of linear solvers.

\begin{figure}
\begin{center}
 \raisebox{0mm}{\includegraphics[trim=10 40 15 20, clip, width=0.45\textwidth]{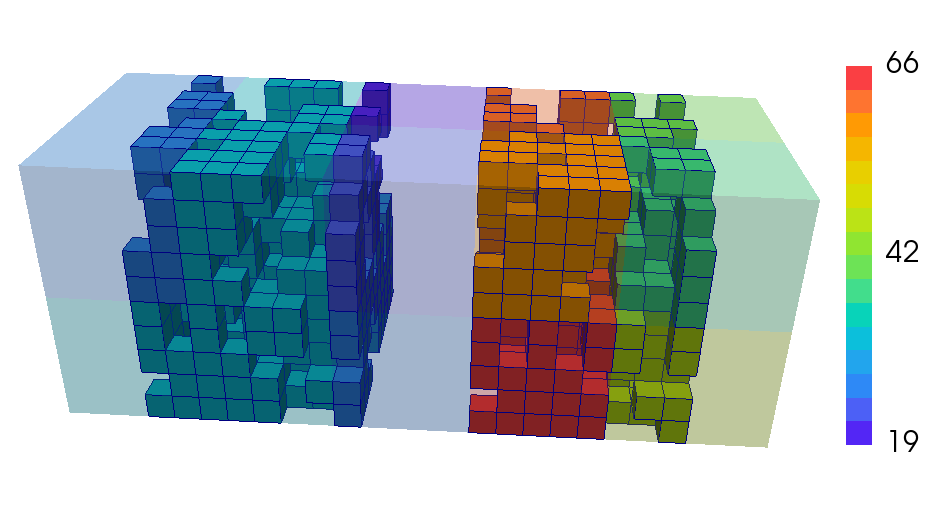}}
 \hfill
 \raisebox{0mm}{\includegraphics[trim=10 40 15 20, clip, width=0.45\textwidth]{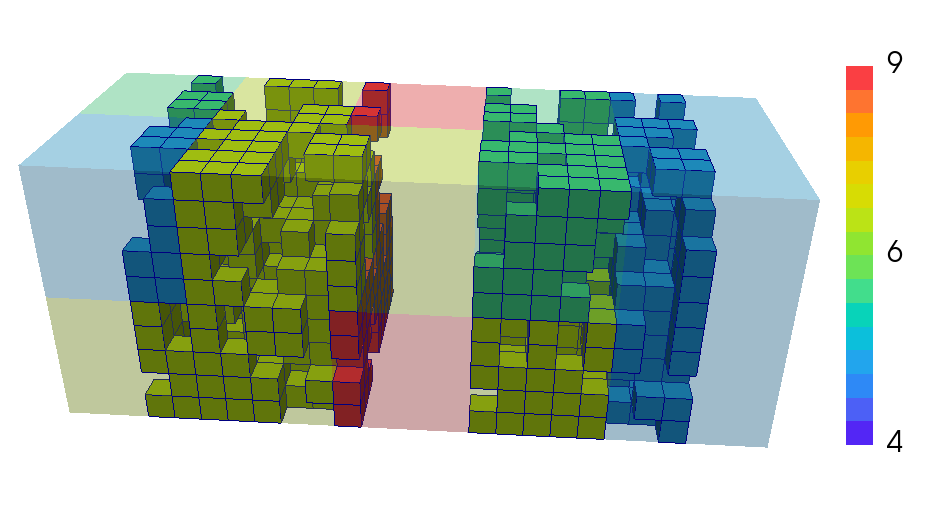}}
\end{center}	
\caption{Maximum local reduced basis dimensions $\dim(\tilde{V}_{c,i})$ (left) and $\dim(\tilde{V}_{\phi,i})$ (right)
obtained in the numerical experiment.\label{fig:rb_dims}}
\end{figure}

In this contribution we investigate if spatially resolved electrochemistry simulations are
in principle amenable to such localization techniques.
For this we partition the computational domain with a cuboid macro grid with elements
$\Omega_1, \ldots \Omega_K$ that are aligned with the microscale voxel grid of the given finite volume
discretization (cf.\ \cref{fig:rb_dims}).
This partition induces a direct sum decomposition of $V_h$:
\begin{equation*}
	V_h = V_{h,1} \oplus \ldots \oplus V_{h,K}, \qquad
	V_{h_i} = \{f \in V_h\,|\, \operatorname{supp}(f) \subseteq \overline{\Omega_i}\}.
\end{equation*}
We now compute local reduced spaces $\tilde{V}_{c,i}$, $\tilde{V}_{\phi,i}$ by first computing
global solution snapshots $c^{(n)}_{\mu_s}$, $\phi^{(n)}_{\mu_s}$ for preselected
parameters $\mu_1, \ldots, \mu_S$ and then
performing local PODs of the $L^2$-orthogonal projections of these snapshots onto the local finite volume spaces
$V_{h,i}$. Hence,
\begin{gather*}
	\tilde{V}_{c,i} \subseteq \operatorname{span}\{P_{V_{h,i}}(c^{(n)}_{\mu_s})\ |\ 1 \leq s
	\leq S,\ 1 \leq n \leq N\},\\
	\tilde{V}_{\phi,i} \subseteq \operatorname{span}\{P_{V_{h,i}}(\phi^{(n)}_{\mu_s})\ |\ 1 \leq s
		\leq S,\ 1 \leq n \leq N\}.
\end{gather*}
Since our high-dimensional model is already given as a non-conforming discretization, we con now
simply obtain a reduced model by solving \cref{eq:reduced} with the reduced space
\begin{equation*}
	\tilde{V} = (\tilde{V}_{c,1} \oplus \ldots \oplus \tilde{V}_{c,K}) \oplus (\tilde{V}_{\phi,1} \oplus \ldots
	\oplus \tilde{V}_{\phi, K}).
\end{equation*}

In order to obtain a fully localized model, localized treatment of the nonlinearities $A^{(bv)}$, 
$A^{(\ooc)}$ is required as well. Not only will most of the speedup during the offline phase be lost when
the interpolation data is computed without localization.
Global interpolation basis vectors will also induce a coupling between all
local approximation spaces $\tilde{V}_{c,i}$,
$\tilde{V}_{\phi,i}$.
Thus the block sparsity structure of the Jacobians appearing in the Newton update problems for
solving \cref{eq:ei_reduced} is lost, strongly deteriorating reduced solution times.
Moreover, the additional reduced degrees of freedom due to localization can exhibit a destabilizing
effect when not accounted for while generating the interpolation spaces:
in the limit when each subdomain $\Omega_i$ corresponds to a single voxel, we have
$\tilde{V} = V_h \oplus V_h$ whereas the images of the interpolated operators are
only $M^{(bv)}/M^{(\ooc)}$-dimensional with $M^{(bv)}/M^{(\ooc)} \ll \dim (V_h \oplus V_h)$.

As a first approach to localized treatment of the nonlinear operators, we proceed similar
to the reduced basis generation.
We first construct local empirically interpolated operators
$I^{(*)}_{i,M^{(*)}_i}\circ\tilde{A}^{(*)}_{i,M^{(*)}_i}\circ R^{(*)}_{i,M^{\prime(*)}_i}$
($* = bv,\ooc$) by applying the \textsc{EI-Greedy} algorithm to the projected evaluations
\begin{equation*}
	\{P_{V_{h,i}}(A^{(*)}([c^{(n)}_{\mu_s}, \phi^{(n)}_{\mu_s}]^T))\ |\ 1\leq s\leq S,\ 1\leq n \leq N\}.
\end{equation*}
We then obtain the localized interpolated operators
\begin{equation*}
	A^{(*)} \approx I^{(*)} \circ \tilde{A}^{(*)} \circ R^{(*)},
\end{equation*}
where
\begin{gather}
	I^{(*)}  = \Bigl[I^{(*)}_{1,M^{(*)}_1}, \ldots, I^{(*)}_{K,M^{(*)}_K}\Bigr], \quad
	\tilde{A}^{(*)} =\operatorname{diag}\Bigl(\tilde{A}^{(*)}_{1,M^{(*)}_1}, \ldots, \tilde{A}^{(*)}_{K,M^{(*)}_K}\Bigr),
	\label{eq:def_I_K} \\
	R^{(*)} = \Bigl[R^{(*)}_{1,M^{\prime(*)}_1}, \ldots, R^{(*)}_{K,M^{\prime(*)}_K}\Bigr]^T. \label{eq:def_R}
\end{gather}
Using these operators in \cref{eq:ei_reduced} leads to a basic, fully localized and fully reduced approximation
scheme for \cref{eq:detailed}.

In order to obtain a stable reduced scheme, accurate approximation of the Butler-Volmer fluxes
is crucial.
However, each localized interpolated operator only takes interface fluxes into its associated
domain $\Omega_i$ into account:
Let $T_1$ be a finite volume cell at the boundary of $\Omega_i$
and $T_2$ an adjacent cell in a different subdomain $\Omega_j$, $i \neq j$.
Unless both cells are selected as interpolation points for the respective operators,
local mass conservation will be violated at the $T_1$/$T_2$ interface due to the
errors introduced by empirical interpolation.

To investigate whether these jumps in the interface fluxes of the interpolated operators
have a destabilizing effect, we consider the following modified scheme:
We denote by $A^{\prime(*)}_i: V_h \oplus V_h \to V_h \oplus V_h$, ($* = bv,\ooc$) the operator
which accumulates all numerical fluxes associated with $A^{(*)}$ which correspond to
grid faces contained in $\overline{\Omega_i}$. Fluxes corresponding to faces which are also
contained in some $\overline{\Omega_j}$, $i \neq j$, are scaled by $1/2$.
This scaling ensures that we have
\begin{equation*}
	A^{(*)} = \sum_{i=1}^K A^{\prime(*)}_i.
\end{equation*}
Each operator $A^{\prime(*)}_i$ is interpolated separately yielding approximations
$I^{\prime(*)}_{i,M^{(*)}_i}\circ\tilde{A}^{\prime(*)}_{i,M^{(*)}_i}\circ R^{\prime(*)}_{i,M^{\prime(*)}_i}$,
where the interpolation data is again obtained via \textsc{EI-Greedy} algorithms for
the evaluations
\begin{equation*}
	\bigl\{A^{\prime(*)}([c^{(n)}_{\mu_s}, \phi^{(n)}_{\mu_s}]^T)\ |\ 1\leq s\leq S,\ 1\leq n \leq N\bigr\}.
\end{equation*}
We then proceed as before by defining $I^{\prime(*)}$, $\tilde{A}^{\prime(*)}$, $R^{\prime(*)}$ as in 
\cref{eq:def_I_K,eq:def_R}, obtaining the localized approximation 
$A^{(*)} \approx I^{\prime(*)} \circ \tilde{A}^{\prime(*)} \circ R^{\prime(*)}$.

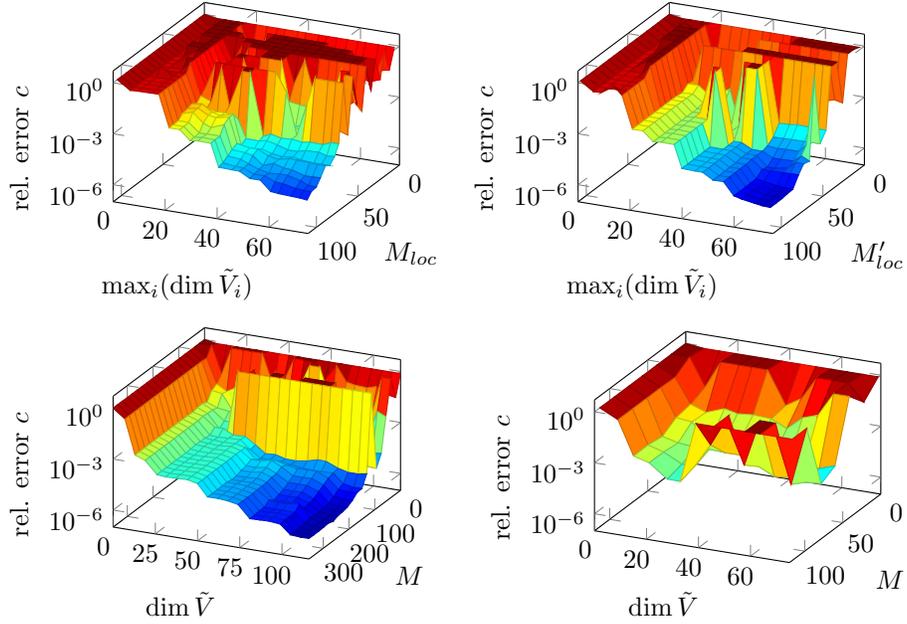
\begin{figure}[t]

	\begin{center}
                \begin{tikzpicture}
			\begin{axis}[width=5.4cm,
				     xlabel={$\max_i(\dim \tilde{V}_i$)},
				     xtick={0, 20, 40, 60},
				     xlabel shift=-5pt,
				     ylabel={$M_{loc}$},
				     ylabel style={yshift=5pt},
				     ylabel shift=-5pt,
				     zlabel={rel. error $c$},
				     ztick={1e0, 1e-3, 1e-6},
				     zlabel style={xshift=4pt},
				     y dir=reverse,
				     zmode=log,
				     colormap/jet,
			             restrict x to domain=0:80,
				     xmax=75,
			             restrict y to domain=0:110,
				     ymax=110,
			     	     zmin=1e-7,
				     point meta min=-14,
				     point meta max=0,
			            ]
				\addplot3[surf] file {lrb_simple_c.dat};
			\end{axis}
		\end{tikzpicture}
		\hfill
		\begin{tikzpicture}
			\begin{axis}[width=5.4cm,
				     xlabel={$\max_i(\dim \tilde{V}_i$)},
				     xlabel shift=-5pt,
				     xtick={0, 20, 40, 60},
				     ylabel={$M^\prime_{loc}$},
				     ylabel style={yshift=5pt},
				     ylabel shift=-5pt,
				     zlabel={rel. error $c$},
				     ztick={1e0, 1e-3, 1e-6},
				     zlabel style={xshift=4pt},
				     y dir=reverse,
				     ztick={1e0, 1e-3, 1e-6},
				     zmode=log,
				     colormap/jet,
			             restrict x to domain=0:80,
				     xmax=75,
			             restrict y to domain=0:110,
				     ymax=110,
			     	     zmin=1e-7,
				     point meta min=-14,
				     point meta max=0,
			            ]
				\addplot3[surf] file {lrb_improved_c.dat};
			\end{axis}
		\end{tikzpicture}
		\\[\medskipamount]
	\begin{tikzpicture}
			\begin{axis}[width=5.4cm,
				     xlabel={$\dim \tilde{V}$},
				     xlabel shift=-5pt,
				     xtick={0, 25, 50, 75, 100},
				     ylabel={$M$},
				     ylabel style={yshift=5pt},
				     ylabel shift=-3pt,
				     ytick={0, 100, 200, 300},
				     zlabel={rel. error $c$},
				     ztick={1e0, 1e-3, 1e-6},
				     zlabel style={xshift=4pt},
				     y dir=reverse,
				     zmode=log,
			             colormap/jet,
			     	     zmin=1e-7,
				     point meta min=-14,
				     point meta max=0,
			            ]
				\addplot3[surf] file {rb_c.dat};
			\end{axis}
		\end{tikzpicture}
		\hfill	
                \begin{tikzpicture}
			\begin{axis}[width=5.4cm,
				     xlabel={$\dim \tilde{V}$},
				     xlabel shift=-5pt,
				     xtick={0, 20, 40, 60},
				     ylabel={$M$},
				     ylabel style={yshift=5pt},
				     ylabel shift=-3pt,
				     zlabel={rel. error $c$},
				     ztick={1e0, 1e-3, 1e-6},
				     zlabel style={xshift=4pt},
				     y dir=reverse,
				     zmode=log,
				     colormap/jet,
			             restrict x to domain=0:75,
				     xmax=75,
			             restrict y to domain=0:120,
				     ymax=120,
			     	     zmin=1e-7,
				     point meta min=-14,
				     point meta max=0,
				    ]
				\addplot3[surf] file {rb_c.dat};
			\end{axis}
		\end{tikzpicture}
	\end{center}
        \caption{%
		Relative model order reduction errors for the concentration variable $c$.
		The error is measured in the $L^2$-in space, $L^\infty$-in time,
		$L^\infty$-in $\mu$ norm for 10 randomly sampled parameters $\mu \in \mathcal{P} := [0.00012, 0.0012]\, A/cm^2 \approx
		[0.1, 1]\, C$. Top left: errors for the fully localized scheme, $\tilde{V}_i:= \tilde{V}_{c,i} \oplus
		\tilde{V}_{\phi,i}$, $M_{loc}:= \max_i (\max(M^{(bv)}_i, M^{(\ooc)}_i))$.
		Top right: errors for the fully localized scheme with additional special treatment of the
		interface fluxes, $M^\prime_{loc}:= \max_i (\max(M^{\prime(bv)}_i, M^{\prime(\ooc)}_i))$.
		Bottom left: errors for reduced basis approximation without localization, $M:=\max(M^{(bv)},
		M^{(\ooc)})$. Bottom right: errors for reduced basis approximation without localization with same
		axis scaling as in top row.
\label{fig:errors}
}
\end{figure}

\section{Numerical experiment}\label{sec:experiment}

As a first numerical experiment we consider the small battery geometry depicted in \cref{fig:geometry}.
For this geometry we compare the performance of the two localized RB approximation schemes introduced 
in \cref{sec:lrb} with the scheme without localization described in \cref{sec:rb}.

The model was simulated for 2000 seconds with equidistant time steps of size $\Delta t:=10s$.
In order to preclude any effects from possibly insufficient sampling of the solution manifold,
the reduced models were constructed using the relatively large amount of $S=20$ equidistant parameters
in the parameter domain $\mathcal{P}:= [0.00012, 0.0012]\, A/cm^2 \approx [0.1, 1]\, C$.
All reduced approximation and interpolation spaces were computed with relative POD/\textsc{EI-Greedy}
error tolerances of $10^{-7}$.
The resulting local reduced basis dimensions for the concentration and potential variables are
depicted in \cref{fig:rb_dims}.
The maximum model reduction errors were estimated by computing the reduction errors for a test
set of 10 random parameters and are shown for the concentration variable in \cref{fig:errors}
(the errors in the potential variable show similar behaviour).
All simulations of the high-dimensional finite volume battery model have been performed within the DUNE software framework
\cite{MR2421580,MR2421579}, which has been integrated with our model order reduction library \texttt{pyMOR} \cite{MRS15}.

We observe (\cref{fig:errors}, top row) that both localized schemes yield stable reduced order models with good error decay, provided
a sufficiently large number of interpolation points is chosen. 
The localized scheme with special treatment of the boundary fluxes (top right) is indeed overall more
stable than the localized scheme without boundary treatment (top left) and yields slightly smaller reduction
errors.

In comparison to the global RB approximation (bottom left), less reduced basis vectors/interpolation points are required
per subdomain to obtain a good approximation for the localized schemes. 
As expected for localized RB schemes, the total number of basis vectors/interpolation points
is larger (cf.\ bottom right) than for the global scheme, however. Given the small size of the
full order model, we cannot expect any speedup for the localized reduced models.
Nevertheless, based on our experience with global RB approximation of this model \cite{OhlbergerRaveEtAl2016},
we expect only a small increase in the number of required basis vectors/interpolation points to approximate
larger, finely resolved geometries. Thus, good speedups can be expected for large-scale applications.
Verifying this hypothesis, as well as developing algorithms for efficient localized construction
and enrichment of the local approximation spaces, will be subject to future work.

\section{Conclusion}
In this contribution we demonstrated the applicability of the Localized Reduced Basis 
Method for an instationary nonlinear finite volume Li-ion battery model with resolved pore scale electrode geometry. 
To this end, we have extended the Localized Reduced Basis Method to parabolic systems of equations, 
while simultaneously employing the localized empirical operator interpolation in order to deal with the strong 
nonlinearities of the underlying electrochemical reaction processes. Numerical experiments were given to 
demonstrate the model order reduction potential of this approach.

\bibliographystyle{spmpsci}
\bibliography{paper}

\begin{thebibliography}{10}
\providecommand{\url}[1]{{#1}}
\providecommand{\urlprefix}{URL }
\expandafter\ifx\csname urlstyle\endcsname\relax
  \providecommand{\doi}[1]{DOI~\discretionary{}{}{}#1}\else
  \providecommand{\doi}{DOI~\discretionary{}{}{}\begingroup
  \urlstyle{rm}\Url}\fi

\bibitem{AlbrechtHaasdonkEtAl2012}
Albrecht, F., Haasdonk, B., Ohlberger, M., Kaulmann, S.: The localized reduced
  basis multiscale method.
\newblock Proceedings of Algoritmy 2012, Conference on Scientific Computing,
  Vysoke Tatry, Podbanske, September 9-14, 2012 pp. 393--403 (2012)

\bibitem{BarraultMadayEtAl2004}
Barrault, M., Maday, Y., Nguyen, N.C., Patera, A.T.: An `empirical
  interpolation' method: application to efficient reduced-basis discretization
  of partial differential equations.
\newblock C. R. Math. Acad. Sci. Paris \textbf{339}(9), 667--672 (2004)

\bibitem{MR2421580}
Bastian, P., Blatt, M., Dedner, A., Engwer, C., Kl{\"o}fkorn, R., Kornhuber,
  R., Ohlberger, M., Sander, O.: A generic grid interface for parallel and
  adaptive scientific computing. {II}. {I}mplementation and tests in {DUNE}.
\newblock Computing \textbf{82}(2-3), 121--138 (2008)

\bibitem{MR2421579}
Bastian, P., Blatt, M., Dedner, A., Engwer, C., Kl{\"o}fkorn, R., Ohlberger,
  M., Sander, O.: A generic grid interface for parallel and adaptive scientific
  computing. {I}. {A}bstract framework.
\newblock Computing \textbf{82}(2-3), 103--119 (2008)

\bibitem{BuhrEngwerEtAl2015}
Buhr, A., Engwer, C., Ohlberger, M., Rave, S.: Application of the {ArbiLoMod}
  to problems of electrodynamics.
\newblock In: Proceedings of MoRePaS III (2015)

\bibitem{BuhrEngwerEtAl2014}
Buhr, A., Engwer, C., Ohlberger, M., Rave, S.: {ArbiLoMod}, a simulation
  technique designed for arbitrary local modifications.
\newblock arXiv e-prints (1512.07840) (2015).
\newblock Http://arxiv.org/abs/1512.07840 (submited to SISC)

\bibitem{DrohmannHaasdonkEtAl2012}
Drohmann, M., Haasdonk, B., Ohlberger, M.: Reduced basis approximation for
  nonlinear parametrized evolution equations based on empirical operator
  interpolation.
\newblock SIAM J. Sci. Comput. \textbf{34}(2), A937--A969 (2012)

\bibitem{HaasdonkOhlbergerEtAl2008}
Haasdonk, B., Ohlberger, M., Rozza, G.: A reduced basis method for evolution
  schemes with parameter-dependent explicit operators.
\newblock Electron. Trans. Numer. Anal. \textbf{32}, 145--161 (2008)

\bibitem{HeinLatz2016}
Hein, S., Latz, A.: Influence of local lithium metal deposition in 3d
  microstructures on local and global behavior of lithium-ion batteries.
\newblock Electrochimica Acta \textbf{201}, 354 -- 365 (2016)

\bibitem{LatzZausch2011}
Latz, A., Zausch, J.: Thermodynamic consistent transport theory of li-ion
  batteries.
\newblock Journal of Power Sources \textbf{196}(6), 3296 -- 3302 (2011)

\bibitem{MRS15}
Milk, R., Rave, S., Schindler, F.: {pyMOR} - generic algorithms and interfaces
  for model order reduction.
\newblock Accepted for publication in SIAM J. Sci. Comput.  (2016)

\bibitem{OhlbergerRaveEtAl2015}
Ohlberger, M., Rave, S., Schindler, F.: True error control for the localized
  reduced basis method for parabolic problems.
\newblock In: Proceedings of MoRePaS III (2015)

\bibitem{OhlbergerRaveEtAl2016}
Ohlberger, M., Rave, S., Schindler, F.: Model reduction for multiscale
  lithium-ion battery simulation.
\newblock In: Numerical Mathematics and Advanced Applications - ENUMATH 2015
  (to appear), LNCSE. Springer (2016).
\newblock DOI 10.1007/978-3-319-39929-4

\bibitem{OhlbergerRaveEtAl2014}
Ohlberger, M., Rave, S., Schmidt, S., Zhang, S.: A model reduction framework
  for efficient simulation of li-ion batteries.
\newblock In: J.~Fuhrmann, M.~Ohlberger, C.~Rohde (eds.) Finite Volumes for
  Complex Applications VII-Elliptic, Parabolic and Hyperbolic Problems,
  \emph{Springer Proceedings in Mathematics \& Statistics}, vol.~78, pp.
  695--702. Springer International Publishing (2014)

\bibitem{OhlbergerSchindler2015}
Ohlberger, M., Schindler, F.: Error control for the localized reduced basis
  multiscale method with adaptive on-line enrichment.
\newblock SIAM J. Sci. Comput. \textbf{37}(6), A2865--A2895 (2015)

\bibitem{PopovVutovEtAl2011}
Popov, P., Vutov, Y., Margenov, S., Iliev, O.: Finite volume discretization of
  equations describing nonlinear diffusion in li-ion batteries.
\newblock In: I.~Dimov, S.~Dimova, N.~Kolkovska (eds.) Numerical Methods and
  Applications, \emph{Lecture Notes in Computer Science}, vol. 6046, pp.
  338--346. Springer Berlin Heidelberg (2011)

\bibitem{sirovich87}
Sirovich, L.: Turbulence and the dynamics of coherent structures. i-coherent
  structures. ii-symmetries and transformations. iii-dynamics and scaling.
\newblock Quarterly of Applied Mathematics \textbf{45}(3), 561--571 (1987)

\end{thebibliography}

\end{document}